\documentclass{amsart}
\usepackage{amssymb, amsmath, latexsym, amscd}
\usepackage{graphics}
\newtheorem{theorem}{Theorem}[section]
\newtheorem{corollary}[theorem]{Corollary}
\newtheorem{lemma}[theorem]{Lemma}

\newtheorem{proposition}[theorem]{Proposition}

\theoremstyle{definition}
\newtheorem{definition}[theorem]{Definition}

\def\bdry{\partial}
\def\bG{\bdry_E}
\def\bE{\bdry_H}
\newcommand{\G}{X}
\newcommand{\GH}{Y}
\newcommand{\GK}{Z}
\newcommand{\GL}{W}
\newcommand{\R}{\mathbb R}
\newcommand{\g}{\mathcal G}
\newcommand{\h}{\mathcal H}

\begin{document} 
\title{Infinitesimal  Operations on  Complexes of Graphs}
\author[J. Conant]{Jim Conant}
\address{Cornell University}
\author[K. Vogtmann]{Karen Vogtmann}
\address{Cornell University}
\thanks{The first author was partially supported by NSF VIGRE grant 
DMS-9983660. The second author was partially supported by NSF grant
DMS-9307313}
\subjclass{17B62, 17B63, 17B70, 20F28, 57M07, 57M15, 57M27 }

\begin{abstract}
In two seminal papers Kontsevich used a construction called 
\emph{graph homology} as a bridge  between certain 
infinite dimensional Lie algebras and various topological objects, including
 moduli spaces of curves,  the group of
outer automorphisms of a free group, and invariants of odd dimensional  manifolds.
In this paper, we  show that Kontsevich's graph complexes, which include graph complexes
studied earlier by Culler and Vogtmann and by Penner,
have a rich algebraic structure.  We define a Lie bracket and cobracket on
graph complexes,
and in fact show that they are Batalin-Vilkovisky
algebras, and therefore Gerstenhaber algebras. We also find 
 natural subcomplexes on which the bracket and cobracket are 
compatible as a Lie bialgebra.

 Kontsevich's graph
complex construction was  generalized to the context of operads by Ginzburg and
Kapranov, with later generalizations by Getzler-Kapranov and Markl. In
\cite{exposition}, we show that Kontsevich's  results in fact extend  to   general
cyclic operads.    For some operads, including the examples
associated to moduli space and outer automorphism groups of free groups, the subcomplex on
which we have a Lie bi-algebra structure is quasi-isomorphic to the entire connected graph complex.  
In the present paper we  show that all of the new
algebraic operations canonically vanish when the homology functor is
applied, and we expect that the resulting constraints will be useful in
studying the homology of the mapping class group, finite type manifold
invariants and the homology of $Out(F_n)$.

\end{abstract}

\maketitle

\section{Introduction}

  In \cite{kontsevich} and \cite{kontsevich2}, M. Kontsevich 
investigated three ``worlds," or {\it  operads}, which he called   {\it
commutative}, {\it associative}, and {\it Lie}.   For
 each of these operads  he defined an infinite dimensional symplectic Lie
algebra and a chain complex of graphs, and then used invariant theory to
prove that the graph complex  computes  the homology of the Lie algebra.
Ginzburg and Kapranov generalized the notion of graph complex to the case of an arbitrary
operad, calling the result the \emph{cobar complex} \cite{ginzkap}.  
Later E. Getzler and M. Kapranov \cite{gk2} generalized Kontsevich's graph
complex construction to the class of  differential graded modular
operads, and called the resulting functor to graph complexes the
\emph{Feynman transform}.  M. Markl 
 \cite{markl}  also  gave a construction  of graph complexes, in the
context of cyclic operads. The purpose of the present paper is to
demonstrate that these graph complexes, i.e. the image of the
Getzler-Kapranov Feynman transform, carry a rich algebraic structure. In
\cite{newconant}, much of this structure is extended to give additional,
higher-order algebraic operations.

For the associative operad,   Kontsevich's graph complex is the same as one defined by Penner
\cite{penner} to study moduli spaces of punctured surfaces, whereas in the Lie case, the graph complex
comes from Culler and Vogtmann's  ``outer space," which they used to study the group of outer automorphisms
of a free group (see \cite{cv}). The commutative operad gives rise to what Kontsevich
refers to as  ``graph homology."
Homology classes which correspond to trivalent graphs   parameterize  finite type 3-manifold 
invariants (see, e.g.,\cite{Arhus}, \cite{KT},\cite{LMO}). 
The homology in other degrees parameterizes invariants of manifolds of higher odd
dimension. 
In the commutative case, the associated Lie algebra $c_\infty$ can be identified with the direct
limit of  Lie algebras $c_n$, where $c_n$ is the  Lie
algebra of polynomial functions on $\mathbb R^{2n}$ with no linear or
constant terms,  under the standard Poisson bracket.  Alternatively,
$c_n$ can be described as the Lie algebra of derivations of a
polynomial  algebra  which preserve the symplectic form. The
equivalence of these two descriptions comes from the
correspondence,   given by the
symplectic form, between  the Lie algebra of functions and the
Lie algebra of vector fields on $\mathbb R^{2n}$.

The commutative graph complex is spanned by oriented graphs,
where the orientation can be most easily described as an equivalence class
of certain labellings of edges and vertices. The chain complex is
graded by the number of vertices in a
graph, and
the boundary operator $\partial_E$ is given by summing over all edge
contractions. The appropriate notion of induced orientation guarantees that the
square of $\partial_E$ is zero.

After examining Kontsevich's paper closely, we discovered the implicit
presence of another boundary operator $\partial_H$, which anticommutes with
$\partial_E$. It showed up as an error term in a certain diagram that
needed to commute, and represented an oversight in Kontsevich's argument.
In \cite{exposition}, we repair the gap in the more general context of
\emph{cyclic operads.}

This boundary operator is defined by contracting over
pairs of \emph{half-edges.} (See Figure~\ref{halfcontract}, which
depicts the contraction
of the two half edges $h$ and $k$.)
 The commutative graph complex is a Hopf algebra, with multiplication
given by disjoint union
    and comultiplication defined as $1\otimes \G +\G\otimes 1$
on connected graphs and extended multiplicatively.
It is easy to see that
$\partial_E$ is both a derivation and a coderivation with respect to these
operations. However, $\partial_H$ is neither of these. Instead, it
satisfies the Batalin-Vilkovisky axiom, which implies that the defect from
being a derivation is a Lie bracket.  Similarly the deviation from being a
coderivation is a Lie cobracket on graphs. We at first expected these
operations to fit together as a Lie bialgebra, but it turns out they are
only compatible on the subcomplex of connected graphs with no separating
edges. In the Lie and associative cases this subcomplex carries the
homology, as we show in \cite{exposition}.
Computer calculations of F. Gerlits \cite{gerlits} indicate that this 
is not so in the case we concentrate on in this paper, the commutative
case. However, in \cite{cv3}  we show that our operations induce a Lie
bialgebra structure on an appropriate quotient complex  which does indeed
carry the homology.

These operations give 
a rich algebraic structure to the functor assigning chain complexes to cyclic operads.
In the world of Lie algebras, there is a similar natural functor sending a Lie
algebra $\mathfrak g$ to the exterior algebra $\Lambda \mathfrak g$,
considered to be a chain complex with the Chevalley-Eilenberg
differential. This is endowed with the Schouten bracket, which is also 
killed by applying the homology functor. In the last section of this
paper we prove that our bracket comes from the Schouten
bracket on the exterior algebra
$\Lambda c_\infty$.

Sullivan and Chas have studied
similar
algebraic structures on the homology of free loop spaces.
(\cite{chas-sullivan},\cite{chas}). In particular they find a Batalin-Vilkovisky structure,
a Lie bialgebra structure, and an uncountable infinity of $Lie_\infty$ structures.
These operations are generically present at the chain complex level but do not form a BV
algebra until passing to homology.

The new algebraic operations described in this paper may be useful in 
studying the three topological applications of Kontsevich's theory:
finite type invariants of odd dimensional manifolds, the
homology of the mapping class group and the homology of $Out(F_n)$.
The vanishing of these operations homologically imposes serious
constraints on these objects, and in particular may be useful in
obtaining information about dimension. 

We also hope that the new operations will be useful in the analysis of
the Feynman transform functor, which Getzler and Kapranov show  
is a homotopy equivalence between certain categories of modular operads.

In this paper we study the properties of the new differential $\partial_H$ for the
commutative graph complex.
 In
\cite{newconant} and
\cite{exposition} we explain the modifications of this paper which are necessary to define
bracket and cobracket for general cyclic operads.

We would like to thank F. Gerlits, E. Getzler, S. Mahajan, D. Sullivan,
and D. Thurston for their interest and stimulating conversations.

\section{Chain complexes of graphs}

In this section we describe Kontsevich's commutative graph complex $\g$ and
the two  boundary operators, $\bG$ and $\bE$.

    By a {\it graph} we mean a finite 1-dimensional $CW$-complex $\G$,
with vertices $v(\G)$ and edges $e(\G)$.  We assume that all vertices
in the graph
have valence at least 3. An {\it orientation} on a graph $\G$ is
simply an orientation of
the vector space
$\mathbb R^{e(\G)}\times H_1(\G;\R)$.  We will usually find it more convenient
to think of an orientation as an equivalence class of labellings, where
 a labelling of $\G$ consists of an ordering of  the vertices $v(\G)$
and arrows on all edges.  Reversing the arrow on any edge, or switching the
order of two vertices changes the  orientation.
These two notions are equivalent for connected graphs (See
\cite{thurston},
\cite{KT} for an explanation of this equivalence).

The
$k$-chains of $\g$ are linear combinations of oriented graphs $(\G,or)$ with
$k$ vertices, modulo the relation $(\G,or)=-(\G,-or)$. This relation
forces all graphs with loops
to be zero, since one can switch the arrow on the loop to get an
isomorphic graph with the opposite orientation, giving
$(\G,or)=-(\G,or)$. Therefore we may assume that our graphs have no loops.

Given an edge $e$ of $\G$, we define $\G_e$ to be the graph obtained
from $\G$ by collapsing $e$ to a point. The first boundary operator
$\bG$ is given by summing over all possible edge collapses:
\begin{definition} Let  $(\G,or)$ be an oriented graph.  Then
     $$\bG(\G, or)=\sum (\G_e,or),$$ where the sum is over all edges 
$e$ of $\G$,
and $\G_e$ is given the
orientation induced from the orientation on $\G$.
\end{definition}
To specify the induced orientation on $\G_e$ in terms of labelled
graphs, choose a representative for the
orientation of
$\G$ such that the initial vertex of
$e$ is labelled 1 and the terminal vertex is labelled 2.  The
labelling on $\G_e$ is then given by the following rule:
the vertex which results from collapsing $e$ is numbered
$1$  and the numbering on all other vertices is reduced by one.  The
arrows on all uncontracted edges  are unchanged.

\begin{lemma} \label{ijcollapse} Choose a labelling to represent the 
oriented graph $\G$.  Then
collapsing an edge of $\G$ from vertex
$i$ to vertex $j$ with $i<j$ induces the orientation $(-1)^jor$, 
where $or$ is the orientation which results from numbering the 
collapsed edge $i$
and reducing the numbering on the vertices labelled $j,\ldots,n$ by 1.
\end{lemma}

Since
all of our graphs will be oriented, we will suppress the orientation
in our notation, writing simply $\G$
instead of $(\G,or)$.

\begin{figure}[ht]
\begin{center}
\includegraphics{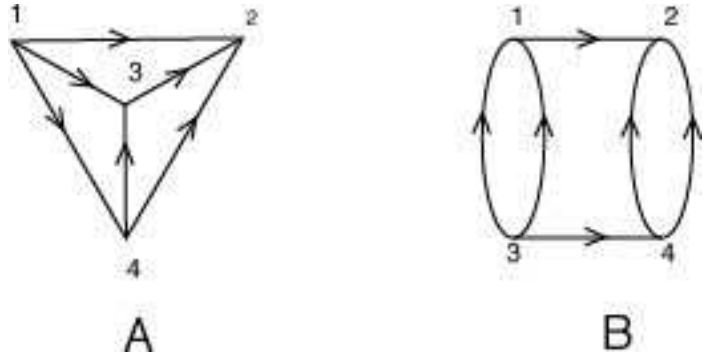}
\caption{Graphs A and B}\label{ABgraphs}
\end{center}
\end{figure}

Since the boundary operator preserves the first Betti number, or
{\it loop degree\,} of a graph, one can decompose graph homology as a
direct sum over the homologies of a fixed loop degree.
For loop degree two, there is
only one possible graph, the theta graph. (All
other possibilities are excluded because they have loops or vertices
of valence 1 or 2.)
Therefore the theta graph gives rise
to a degree two homology class. In loop degree three, there are two
possible graphs which have four vertices, as shown in Figure \ref{ABgraphs}.

\begin{figure}[ht]
\begin{center}
\includegraphics{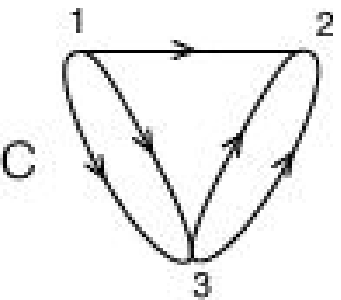}
\caption{Graph C}\label{Cgraph}
\end{center}
\end{figure}

There is only one graph with three vertices and no loops, shown in 
Figure \ref{Cgraph}

There is also a graph on two vertices with no loops, but this has an 
orientation
reversing automorphism. It is easy to see
that $\partial_E A = -6C$ and
$\partial_E B = 2C$. Hence we get one homology class in degree four,
$A+2B$.

\begin{figure}[ht]\begin{center}
\includegraphics{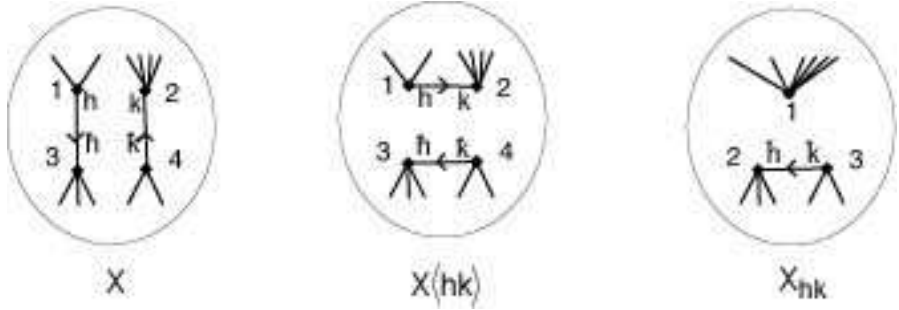}
\caption{Contracting half-edges $h$ and $k$}\label{halfcontract}
\end{center}
\end{figure}

To describe the second boundary operator on $\g$, we use the  {\it
half-edges} of a graph $\G$.  Each half-edge $h$ begins at a vertex
$v(h)$, is contained in an edge
$e(h)$ and  has a complementary half-edge
$\bar h$, with $h\cup\bar h =e(h)$.   Given two half-edges $h$ and
$k$   of $\G$, we form a new graph
$\G\langle hk \rangle$ by cutting and pasting, as follows:   if 
$k=\bar h$, then
$\G\langle hk \rangle=\G$; if $k\neq \bar h,$ we  cut  to separate
$h$ from $\bar h$ and $k$ from $\bar k$, then glue $h$ to $k$  and
$\bar h$ to $\bar k$ to form two new edges (see Figure~\ref{halfcontract}).
In terms of labelled graphs, the orientation on  $\G\langle hk 
\rangle$ is given
as follows: choose a representative for the orientation on
$\G$ so that
$h$ is the initial half-edge of   $e(h)$, and $k$ is the terminal
half-edge of $e(k)$.   In $\G\langle hk \rangle$, the edge
$h\cup k$ is oriented from $h$ to $k$, and the edge $\bar h\cup\bar
k$ is oriented from $\bar k$ to $\bar h$.

If $h$ and $k$ are half-edges of $\G$ with $v(h)\neq v(k)$, then
$h\cup k$ forms an edge of $\G\langle hk \rangle$, which we can now 
collapse; the
result, $(\G\langle hk \rangle)_{h\cup k}$, is more simply denoted
$\G_{hk}$.   If, on the other, hand, $v(h)=v(k)$, then $\G\langle 
hk\rangle$ has a loop, so is equal to $0$; thus we define
$\G_{hk}$ to be $0$.  Note that
$\G_{hk}=\G_{kh}$ as oriented graphs.

     The second boundary map  on
$\g$ is   given by

\begin{definition}
$$\bE \G=\sum\G_{hk},$$
where the sum is over all pairs $\{h,k\}$ of  half-edges of $\G$ with
$h\neq \bar k$, and $\G_{hk}$ is given the orientation induced from
$\G$.
\end{definition}

To check that $\bG$ and $\bE$ are boundary operators, we use the
following lemma:

\begin{lemma}\label{orient} (Orientation Lemma)
For any four distinct half-edges $h, k, r, s$ of $\G$,
      $$(\G_{hk})_{rs}=-(\G_{rs})_{hk}.$$
\end{lemma}

\begin{proof}

$(\G_{hk})_{rs}$ and $(\G_{rs})_{hk}$ are  the results
of collapsing the edges $h\cup k$ and $r\cup s$ of
$\G\langle hk\rangle\langle rs\rangle =\G\langle rs\rangle \langle hk\rangle $,
in the opposite order.   Now observe that collapsing two edges of an
oriented graph in  opposite order results in isomorphic
graphs with opposite orientations, using Lemma \ref{ijcollapse}.  \end{proof}

\begin{proposition} $\bG^2=\bE^2=(\bG+\bE)^2=0$
\end{proposition}

\begin{proof}

We have
     $$(\bG+\bE)^2(\G)=\sum  (\G_{hk})_{rs}, $$
where the sum is over all sets $\{h,k\}$ and $\{r,s\}$ of half-edges
with $h, k, r$ and $s$ distinct.
These terms  cancel in pairs, by  Lemma \ref{orient}.

Similarly, the terms in
both  squares $(\bG)^2(\G)$ and $(\bE)^2(\G)$ cancel in pairs:  for
$\bG$, the sum is over
all sets $\{h,\bar h\}$ and $\{k,\bar k\}$ with $\{h,\bar
h\}\neq\{k,\bar k\}$, and for $\bE$ the sum is over
all sets of pairs $\{h,k\}$ and $\{r,s\}$ with $k\neq \bar h$, $s\neq
\bar r$ and $\{r,s\}\neq \{\bar h, \bar k\}$ (equivalently,
$\{h,k\}\neq\{\bar r,\bar s\}$).\end{proof}

\begin{corollary} $\bG\bE=-\bE\bG$
\end{corollary}

\begin{figure}[ht]\begin{center}
\includegraphics{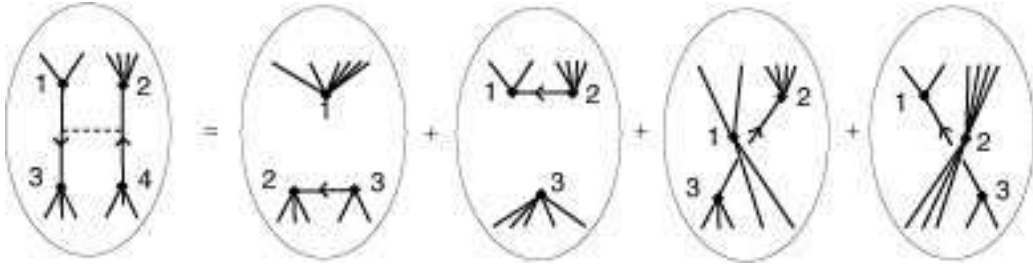} 
\caption{Dotted line notation}\label{dottedline}
\end{center}
\end{figure}

\begin{figure}[ht]
\begin{center}
\includegraphics{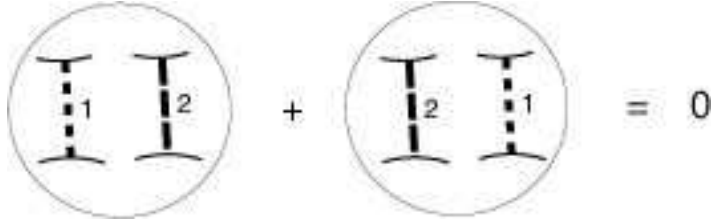}
\caption{Antisymmetry}\label{antisymmetry}
\end{center}
\end{figure}

\begin{figure}[ht]
\begin{center}
\includegraphics{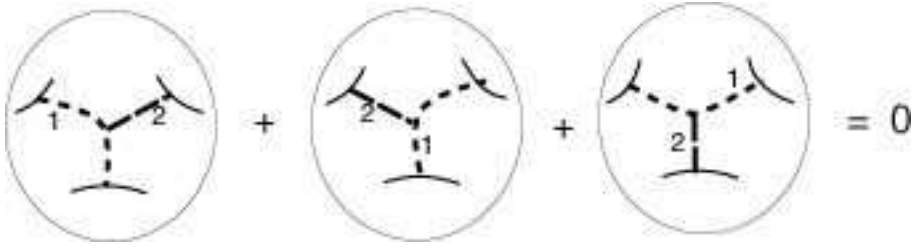}
\caption{Jacobi relation}\label{jacobi}
\end{center}
\end{figure}

We  also briefly mention a slightly different, suggestive visualization
of
$\bE$.
The terms of $\bE \G$ naturally group themselves into sets
of four, namely the four graphs $\G_{hk}$ which can be
formed from the half edges contained in a given pair of edges.
We represent each such set of four graphs graphically by drawing a
dotted line between the corresponding full edges (see Figure \ref{dottedline}).
Now $\bE \G$ is given by summing over all possible ways
of drawing a dotted line between two edges of $\G$. The fact that $\bE^2=0$
can be derived from the   two identities pictured in Figures 
\ref{antisymmetry} and
\ref{jacobi}.    Figure \ref{antisymmetry} represents an antisymmetry 
relation, where the
dotted lines   are numbered to represent the order in
which the operations are performed.
  Figure \ref{jacobi} represents   a sort of Jacobi identity, where a 
second dotted line
coming into a dotted line means that the second dotted line will 
attach to the uncontracted
edge in each summand coming from the first dotted line. When we later define a
    bracket on graphs, these two identities
can be used to give an alternative proof that this bracket
is a Lie bracket.

    Kontsevich's
``associative" and ``Lie" complexes  have similar descriptions
to the above
(``commutative") complex, except that the graphs have additional
structure.   In the associative case graphs come with a cyclic
ordering  of the edges
incident to each vertex.  In the Lie case an equivalence class of 
trivalent trees with $r$
leaves (modulo antisymmetry and  IHX, or ``Jacobi," relations) is 
associated to each vertex of valence
$r$. The operations of cutting, pasting and collapsing described
above can be done in such a way as to
induce natural cyclic orderings or trees on the new vertices created,
so that the boundary operators $\bE$ and
$\bG$ have natural definitions in these settings as well.  In this
paper, we present only the commutative case, for simplicity;
details of the remaining cases will appear in \cite{exposition}.

\section{Some graded algebra \label{graded}}

Both chain complexes and their homology have the structure of graded
vector spaces.   We describe
algebraic structures in this paper which are cognizant of the
grading on the chain complexes of graphs, and descend to
structures on homology. In this section we collect some standard definitions
from graded algebra. We also record how things change under a grading shift.

Suppose $V$ is a graded vector space, either ${\bf Z}$-graded
($V=\oplus_{n\geq 0} V_n$) or ${\bf Z/2}$-graded ($V=V_0\oplus V_1$).
We denote by  $|v|$ the degree of $v\in V$.   A linear map $T\colon
V\to
W$ of graded vector spaces is said to have {\it degree $d$} if
$|T(v)|=|v|+d$ for all homogeneous $v\in V$.

In order to make
algebraic structures derived from
graded vector spaces reflect the grading, we use the {\it twist map}
$\tau\colon V\otimes W\to W\otimes V$, which
takes
$v\otimes w$ to
$(-1)^{|v||w|}w\otimes v$.
This extends to an
isomorphism
$$  V_1\otimes V_2\cdots\otimes V_k\to V_{\pi(1)}\otimes
V_{\pi(2)}\otimes\cdots\otimes
V_{\pi(k)}$$
for any permutation $\pi\in \Sigma_k$ and graded vector spaces
$V_1,V_2,\ldots, V_k$ ; for example,
      the cyclic permutation (123) in $\Sigma_3$ gives an isomorphism
$\sigma\colon V\otimes V\otimes V\to V\otimes V\otimes V$ defined by
$$ \sigma(u\otimes v\otimes w)=(-1)^{|w|(|u|+|v|)} w\otimes u\otimes v.$$
In general, the sign is determined by ``Koszul rule of signs": every
time one switches two adjacent terms in a tensor product, the sign
changes by the product of their degrees.

The isomorphisms given by the Koszul rule   can be interpreted as new
actions of $\Sigma_k$ on the tensor product of $k$ copies of a graded
vector space
$V$: there is a {\it symmetric} action, $$\pi\cdot(v_1\otimes
\cdots\otimes v_k)=(Koszul\  sign)\  v_{\pi(1)}\otimes
\cdots\otimes v_{\pi(k)}$$  and an {\it alternating} action,
$$\pi\cdot(v_1\otimes \cdots\otimes v_k)=(Koszul\  sign)(sign(\pi))\
v_{\pi(1)}\otimes \cdots\otimes
v_{\pi(k)}$$

\begin{definition} The {\it graded wedge product}, denoted
$\Lambda^kV$, is the quotient of $\otimes^k V$ by the alternating
action of $\Sigma_n$,
and the {\it graded symmetric product}, denoted $\odot^kV$ or $S^kV$, is the
quotient of $\otimes^k V$ by the symmetric action of $\Sigma_n$.
\footnote{In the literature $\Lambda V$ often means what we denote by $SV$.}
\end{definition}

\medskip
Let $V[-n]$ denote $V$ with the grading shifted downward by $n$.
The
wedge product and symmetric product are related as follows.

\medskip
\begin{proposition}\label{phi}  There is a natural additive isomorphism $\phi\colon
S^*V\to \Lambda^*\left(V[-1]\right)$, defined inductively by $$\phi(A\odot
B)=(-1)^{|A|}\phi(A)\wedge\phi(B),$$
where $|v_1\odot\ldots\odot v_k|=|v_1|+\ldots+|v_k|$ and
$\phi(v)=v$ for $v\in V$.

In particular, for $k=2$ this isomorphism is given by
      $v\odot w\mapsto (-1)^{|v|}v\wedge w$, and for $k=3$ by $u\odot
v\odot w\mapsto (-1)^{|v|} u\wedge v\wedge w$.
\end{proposition}

The graded definition of Lie algebra \cite{Milnor-Moore} is as follows.
\begin{definition} A {\it graded   Lie bracket} is a linear map
$b\colon V\otimes V\to V$, written $b(x\otimes y)=[x,y],$ satisfying
\begin{itemize}
\item{1.} {\it  (Graded antisymmetry)}
$$[v,w]=-(-1)^{|v||w|}[w,v]$$

\noindent and
\item{2.} {\it  (Graded Jacobi identity)}  $$[u,[v,w]] +
(-1)^{|w|(|u|+|v|)}[w,[u,v]] + (-1)^{|u|(|v|+|w|)}[v,[w,u]]=0$$
\end{itemize}
\end{definition}

\medskip\noindent An equivalent way to state the antisymmetry and
Jacobi relations is that the following two compositions are
zero:
\begin{itemize}
\item{1.} $V\otimes V\ {\buildrel{id +
\tau}\over\longrightarrow}\ V\otimes V\
{\buildrel{b}\over\longrightarrow} \ V$
\item{2.} $ V\otimes V\otimes V\ {\buildrel{id + \sigma +
\sigma^2}\over\longrightarrow}\
      V\otimes V\otimes V\ {\buildrel{id \otimes
b}\over\longrightarrow}\ V\otimes V \ {\buildrel{b
}\over\longrightarrow}\ V.$
\end{itemize}

\medskip
This diagrammatic description is convenient because it allows us to
define a {\it graded  Lie cobracket} by simply reversing all the
arrows:
\begin{definition} A graded Lie cobracket is a linear map
$\theta\colon V\to V\otimes V$ satisfying
\item {1.} {\it  (Graded co-antisymmetry)} $(id + \tau)\circ\theta=0$
\item {2.} {\it  (Graded co-Jacobi identity)}
$(id+\sigma+\sigma^2)\circ(id\otimes\theta)\circ\theta=0$
\end{definition}

Notice that when $V$ is finite-dimensional, $(V,b)$ is a Lie
algebra if and only if
$(V^*,b^*)$ is a Lie coalgebra.  If $V=\oplus V_i$ is a direct sum of 
finite-dimensional
vector spaces, then  $(\oplus V_i,\oplus b_i)$ is a Lie
algebra if and only if
$(\oplus V_i^*,\oplus b_i^*)$ is a Lie coalgebra.   We write 
$V^\dagger=\oplus V_i^*$ and $b^\dagger=\oplus b_i^*$
to avoid confusion with $V^*$ and $b^*$.

\medskip
The  definitions of bracket and cobracket can be reformulated in the
following nice way.
Notice that the antisymmetry conditions imply that the bracket
and cobracket induce maps $$b : V\wedge V\to V$$ and
$$\theta : V\to V\wedge V$$
(In the case of the bracket, we are thinking of $V\wedge V$ as a
quotient of $V\otimes V$, and in the cobracket case as a submodule of
$V\otimes
V$.)
     Further, these maps can be extended
to the entire algebra $\Lambda^* V$ to itself. The map $b$ extending
the bracket is a coderivation, and is just the usual
Lie-algebra-homology boundary map.
The map $\theta$ extends as a derivation.
\begin{lemma}The Jacobi
identity is equivalent to the assertion that $b^2=0$.   The co-Jacobi
identity is
equivalent to the assertion that  $\theta^2=0$ .
\end{lemma}

By the lemma, the Jacobi identity is precisely what
is needed to make the Chevalley-Eilenberg complex $\Lambda^*V$ of a
graded Lie algebra $V$ into a chain
complex.
There is another standard operation on $\Lambda^*V$, the {\it Schouten
bracket}, which usually appears in the context
of Lie algebras of vector fields.

\begin{definition}
Let $V$ be a graded Lie algebra. Then the
    Schouten bracket on $\Lambda^*V$ is defined as follows:

\begin{align*}
&[v_1\wedge\ldots\wedge v_p, w_1\wedge\ldots\wedge w_q]  =\\
     &\qquad\sum_{i,j}(-1)^\epsilon[v_i,w_j]\wedge v_1\wedge\ldots\wedge
\hat{v}_i\wedge\ldots\wedge v_p\wedge w_1\wedge\ldots\wedge
\hat{w}_j\wedge\ldots\wedge w_q,
\end{align*}
where
\begin{align*}
\epsilon  &=  i+j+p+1 + Koszul\\
              &=  i+j+p+1+ |v_i|(|v_1|+\ldots |v_{i-1}|)  \\
             &\qquad\qquad + |w_j|(|v_1|+\ldots\hat{|v_i|}+ \ldots
|v_p|+|w_1|+\ldots |w_{j-1}|)
\end{align*}
\end{definition}
Notice that $i+j+p+1$ is the sign of the permutation bringing
$v_i$,$w_j$ to the
front of the wedge product.

To express the compatibility of bracket and cobracket in the graded
setting we first
review the condition in the ungraded setting. One way of doing this
is as follows
(\cite{majid}),
$$\theta [v,w] = ad_v(\theta (w)) - ad_w (\theta (v)),$$ where the
adjoint action
is extended to the tensor product as a derivation: $ad_v(w_1\otimes w_2)
= [v,w_1]\otimes
w_2 + w_1\otimes[v,w_2]$. This is the same as the following condition:
$$ \theta[v,w] = (b\otimes id +  (id\otimes 
b)\circ\tau_{12})(id\otimes\theta)(id-\tau)(v\otimes w),$$
where
$\tau_{12}$ is the transposition swapping the first two tensor factors. Let
$\theta(v)=\sum v_1\otimes v_2$ and $\theta(w)=\sum w_1\otimes w_2$.
Adding the Koszul signs and the degrees $|b|$ of the bracket and 
$|\theta|$ of the cobracket in the graded situation,
  the above
condition becomes
\begin{align*}
\theta[v,w]&=(b\otimes id +  (id\otimes b)\tau_{12})((-1)^{|\theta||v|}v\otimes w_1\otimes w_2\\
&\qquad\qquad\qquad\qquad\qquad\qquad\qquad\qquad -
(-1)^{|v||w|+|\theta||w|}w\otimes v_1\otimes v_2)\\
&=(-1)^{|\theta||v|}[v,w_1]\otimes
w_2-(-1)^{|v||w|+|\theta||w|}[w,v_1]\otimes v_2
\\  &\quad+ (-1)^{|\theta||v|+|v||w_1|+|b||w_1|} w_1\otimes [v,w_2]
\\ &\qquad-(-1)^{|v||w|+|\theta||w|+|w||v_1|+|b||v_1|}v_1\otimes [w,v_2]
\end{align*}
Now, passing to wedge products,
\begin{align*}
\theta[v,w] &=(-1)^{|\theta||v|}[v,w_1]\wedge
w_2-(-1)^{|v||w|+|\theta||w|}[w,v_1]\wedge v_2 \\
&\quad + (-1)^{|\theta||v|+|v||w_1|+|b||w_1|} w_1\wedge [v,w_2]\\
&\qquad -(-1)^{|v||w|+|\theta||w|+|w||v_1|+|b||v_1|}v_1\wedge [w,v_2]\\
&=(-1)^{|\theta||v|}[v,w_1]\wedge w_2
-(-1)^{|v||w|+|\theta||w|}[w,v_1]\wedge v_2 \\
&\quad +(-1)^{|\theta||v|+|v||w_1|+|b||w_1| +|w_1|(|v|+|w_2|+|b|)+1}[v,w_2]\wedge
w_1\\
&\qquad
-(-1)^{|v||w|+|\theta||w|+|w||v_1|+|b||v_1|+|v_1|(|w|+|v_2|+|b|)+1}
[w,v_2]\wedge v_1\\
&= (-1)^{|\theta||v|}\big([v,w_1]\wedge w_2 +
(-1)^{|w_1||w_2|+1}[v,w_2]\wedge w_1\big)\\
&\quad +(-1)^{|v||w|+|\theta||w|}\big((-1)^{|v_1||w|}[v_1,w]\wedge
v_2 \\
&\qquad+ (-1)^{|v_1||v_2|+|v_2||w|+1}[v_2,w]\wedge v_1\big)\\
&= (-1)^{|\theta||v|}[v,w_1\wedge w_2] + (-1)^{|v_2||w|}[v_1,w]\wedge
v_2\\
&\quad + (-1)^{|v_1||v_2|+|v_1||w|+1}[v_2,w]\wedge v_1\\
&= (-1)^{|\theta||v|}[v,w_1\wedge w_2] -[v_1\wedge v_2,w]
\end{align*}

We therefore adopt the following definition of graded Lie bialgebra.
\begin{definition}
A graded Lie bialgebra is a graded vector space $V$ together
with a Lie bracket $b=[\cdot,\cdot] : \Lambda^2V\to V$ and a Lie
cobracket $\theta : V\to\Lambda^2V$ which are compatible:
$$\theta[v,w] =
-[\theta(v),w] +(-1)^{|v||\theta|}[v,\theta (w)],$$
where the bracket on the right hand side of the equation is the
Schouten bracket.
\end{definition}


We would like to analyze now what happen when the degree is shifted. Recall
(Proposition~\ref{phi}) that there is a natural isomorphism $\phi\colon \Lambda V[-1] \cong S V$. If one has
an operation
$b = [\cdot,\cdot]\colon S^2V\to V$, one gets an induced operation
$b^\phi = [\cdot,\cdot]^\phi\colon \Lambda^2V[-1]\to V$ defined by $b^\phi = b\circ \phi^{-1}$. Similarly,
an operation $\theta\colon V\to S^2V$ gives rise to a map $\theta^\phi=\phi\circ \theta$.

\begin{proposition}\label{prop1} 
A linear map $b=[\cdot ,\cdot]\colon S^2V\to V$ induces a graded Lie bracket
$[\cdot,\cdot]^\phi$ on
$V[-1]$ iff either of the following two equivalent conditions holds:
\begin{itemize}\item{1.} The extension of $b$ to $SV$ as a coderivation has trivial square.
\item{2.} {\it  (Graded Jacobi identity)}  $$[u,[v,w]] +
(-1)^{|w|(|u|+|v|)}[w,[u,v]] + (-1)^{|u|(|v|+|w|)}[v,[w,u]]=0$$
\end{itemize}
\end{proposition}

Note that the second axiom is the usual graded Jacobi identity!

\begin{proposition}
A linear map $\theta\colon V\to S^2V$ induces a graded Lie cobracket $\theta^\phi$ on
$V[-1]$ iff the extension of $\theta$ to $SV$ as a derivation has trivial square.
\end{proposition}

As before, if $V=\oplus V_i$ is a direct sum of finite dimensional vector spaces we have
that $\theta\colon V\to S^2V$ satisfies the graded co-Jacobi identity iff
$\theta^\dagger$ satisfies the graded Jacobi identity. In othr words, $\theta^\phi$ is a Lie
cobracket iff $(\theta^\dagger)^\phi$ is a Lie bracket.

In order to define the compatibility conditions between bracket and cobracket
in the symmetric world, we need to
transfer the Schouten bracket to the symmetric setting.

\begin{definition}
Let $V$ be a graded Lie algebra. The symmetric Schouten
bracket on $S V$ is defined by:

\begin{align*}
&[v_1\odot\ldots\odot v_p, w_1\odot\ldots\odot w_q] =\\
&\qquad      \sum_{i,j}(-1)^\epsilon[v_i,w_j]\odot v_1\odot\ldots\odot
\hat{v}_i\odot\ldots\odot v_p\odot w_1\odot\ldots\odot
\hat{w}_c\odot\ldots\odot w_q
\end{align*}
where
\begin{align*}
\epsilon& =
Koszul\\
&=|v_i|(|v_1|+\ldots |v_{i-1}|) + |w_j|(|v_1|+\ldots
\hat{|v_i|}+\ldots |v_p|+|w_1|+\ldots |w_{j-1}|).
\end{align*}
\end{definition}

\begin{proposition}
The linear maps $b\colon S^2V\to V$ and $\theta\colon V\to S^2V$ induce a graded Lie
bialgebra structure on $V[-1]$ iff the following three conditions hold:
\begin{itemize}
\item[1.] $b^\phi$ is a Lie bracket.
\item[2.] $\theta^\phi$ is a Lie cobracket.
\item[3.] $\theta([v,w]) = -[\theta(v),w] -(-1)^{|v||\theta|}[v,\theta (w)]$
\end{itemize}
\end{proposition}

\section{Product and bracket}

Let $\G$ be a labelled graph with  vertices numbered $1,2,\ldots,x $,
and $\GH$ a labelled graph with  vertices numbered $1,2,\ldots,y$.
Define the product
$\G  \cdot \GH$ to be the  disjoint union of $\G$ and $\GH$, with  the
numbering on vertices of
$\GH$  shifted  by adding $x$ to each, thus becoming $x+1,\ldots, x+y$.
Then we
have
$$\G\cdot \GH=(-1)^{xy}{\GH\cdot \G}.$$
This product extends bilinearly to linear combinations of graphs, turning
    $\g$ into a graded commutative algebra. One may allow
    the empty graph as a basis element of $\g$, since it acts as a unit under
the  disjoint union operation.

\begin{lemma} With respect to this product the   boundary operator
$\bG$ is a graded derivation:
$$\bG(\G\cdot \GH)=\bG(\G)\cdot \GH + (-1)^{x} \G\cdot\bG(\GH).$$
\end{lemma}
\begin{proof} This follows since each term $(\G\cdot 
\GH)_{e}$ of $\bG(\G\cdot \GH)$ is
obtained by collapsing an edge $e$, which is either in $\G$ or is in 
$\GH$.  The sign comes from the fact that if $e$ is an edge of  $\GH$,
then $(\G\GH)_e=(-1)^{x}\G\GH_e$. \end{proof}
\medskip

The second boundary operator $\bE$, on the other hand, is not a
derivation; if $h
\subset   \G$ and $k\subset \GH$, then the term  $(\G\cdot\GH)_{hk}$ 
of   $\bE(\G\cdot
\GH)$ is not in $\bE(\G)\cdot
\GH + (-1)^{x} \G\cdot\bE(\GH)$ .
We define the bracket $[\G,\GH]$  so that it measures how far $\bE$ is
from being a derivation:

\begin{definition}
$[\G,\GH]=\bE(\G\cdot \GH)-\bE(\G)\cdot \GH -(-1)^{x}  \G\cdot\bE(\GH).$
\end{definition}

In other words, the bracket of $\G$ and $\GH$ is the sum of all
graphs obtained by contracting a half-edge of $\G$ with a half-edge of $\GH$:

$$[\G,\GH]=\sum_{h\in \G, k\in \GH} (\G\cdot \GH)_{hk}$$

The bracket obeys symmetry and Jacobi relations as given in the
following two lemmas:

\begin{lemma} Let $\G$ and $\GH$ be graphs with  $x$ and $y$
vertices, respectively.  Then
$[\G,\GH]=(-1)^{xy}[\GH,\G].$
\end{lemma}

\begin{proof}
\begin{align*}
[\G,\GH]&=\sum_{h\in \G, k\in \GH}(\G\cdot \GH)_{hk}\\
&=(-1)^{xy}\sum_{h\in \G, k\in \GH}(\GH\cdot \G)_{hk}\\
&=(-1)^{xy}[\GH,\G].
\end{align*}
\end{proof}

\medskip

\begin{lemma} Let $\G,\GH$ and $\GK$ be graphs with  $x, y,$ and $z$
vertices, respectively.
Then
$$[[\G,\GH],\GK]+(-1)^{z(x +y)}[[\GK,\G],\GH]+(-1)^{y(z+x)}[[\GH,\GK],\G]=0.$$
\end{lemma}

\begin{proof}

We have
\begin{align*}
[[\G,\GH],\GK]&=\sum ((\G\cdot  \GH )_{hk}\cdot \GK)_{rs} \\
&=\sum   ((\G\cdot  \GH\cdot \GK)_{hk})_{rs}
\end{align*}
where $h\in \G,
k\in \GH, s\in \GK$ and
$r\in \G\cup \GH, r\not\in\{h,k\}$.

If $r\in \GH$, then  by   Lemma \ref{orient} $((\G\cdot \GH\cdot 
\GK)_{hk})_{rs}$ cancels with the term
$$((\G\cdot \GH\cdot \GK)_{rs})_{hk}=(-1)^{z(x+y) }((\G\cdot \GH\cdot 
\GK)_{rs})_{hk},$$
which is in
$(-1)^{z(x+y)  }[[\GH,\GK],\G]$.

If
$r\in \G$ then  $((\G\cdot \GH\cdot \GK)_{hk})_{rs}$
cancels with the term $((\G\cdot \GH\cdot \GK)_{rs})_{hk}$ of 
$(-1)^{x(y+z)}[[\GK,\G],\GH]$.

The remaining terms of
$(-1)^{x(y +z)}[[\GK,\G],\GH]$ and $(-1)^{z(x+y)}[[\GH,\GK],\G]$ similarly
cancel in pairs.

\end{proof}

\medskip

The above two lemmas combine with Proposition~\ref{prop1} to give

\begin{proposition} $\g [-1]$ is a graded Lie algebra with bracket $[\cdot,\cdot]^\phi$.
\end{proposition}

\noindent{\bf Remark}.
In terms of the dotted line notation we introduced after the definition of
$\bE$, the bracket $[\G,\GH]$ is   the sum over all possible ways of
drawing a dotted line between an edge of $\G$ and an edge of $\GH$. An
alternate proof of the fact that
the bracket satisfies the Jacobi identity can be given using the antisymmetry
and Jacobi identities of dotted lines.
\medskip

In fact, the bracket gives a stronger structure on
$\g[-1]$. Recall that a
{\it Gerstenhaber algebra\/} (or \emph{graded Poisson algebra}) is a graded commutative,
associative algebra $V$ with a degree -1 Lie bracket,
satisfying
\begin{align*}
 [u,vw] &= [u,v]w + (-1)^{|u|(|v|+1)}v[u,w]
\end{align*}

\begin{proposition}
     $\g[-1]$ is a Gerstenhaber algebra, under the product and bracket defined above.
\end{proposition}
\begin{proof} Let $\G$, $\GH$ and $\GK$ be graphs with 
$x$, $y$ and $z$
vertices, respectively. Recall that the bracket on the shifted graph complex is denoted
$[\cdot,\cdot]^\phi$. 
 Then $[A,B]^\phi = (-1)^a [A,B]$, where $a$ is the number of vertices of A.
Now
\begin{align*}
[\G,\GH\GK]&=\sum_{h\subset \G, k\subset  \GH \GK} ((\G\cdot 
(\GH\cdot \GK))_{hk}\\
          &=\sum_{h\subset \G, k\subset \GH} (\G\cdot \GH \cdot \GK)_{hk} +
                \sum_{h\subset G, k\subset \GK} (\G\cdot \GH \cdot \GK)_{hk}\\
&= \sum ( \G\cdot  \GH)_{hk} \cdot \GK + (-1)^{zy}\sum (\G\cdot 
\GK)_{hk}\cdot \GH\\
&= [\G,\GH] \GK + (-1)^{zy}[\G,\GK]\GH\\
&= [\G,\GH] \GK + (-1)^{y(x-1)}\GH[\G,\GK]
\end{align*}
Multiplying through by $(-1)^x$, and noticing that $|x|=x+1$ and $|y|= y+1$, we see that
$[X,YZ]^\phi = [X,Y]^\phi Z +
(-1)^{|x|(|y|+1)} +  Y[X,Z]^\phi$, as desired. 
\end{proof}


An algebraic structure that has recently gained attention is
a {\it Batalin-Vilkovisky algebra\/} (\cite{getzler}, \cite{chas-sullivan}).
It is a graded commutative algebra together with
with a degree 1
map $\Delta$ satisfying $\Delta^2=0$ and such that
$(-1)^{|v|}\Delta(vw) -
(-1)^{|v|}(\Delta v)w - v\Delta w$ is a
derivation of each variable. That is,
\begin{align*}\Delta(uvw) &= \Delta (uv) w + (-1)^{|u|}u\Delta (vw) +
(-1)^{(|u|-1)|v|}v\Delta(u)\\
&\qquad -(\Delta u) vw - (-1)^{|u|}u(\Delta v) w - (-1)^{|u|+|v|}uv(\Delta w).
\end{align*}
If we consider $\g$ as a ${\bf
Z}/2$-graded vector space, then
$\bE$, which is a degree -1 operator, becomes ``degree 1,"  and makes
$\g$ into a super Batalin-Vilkovisky algebra.

\begin{proposition}
As a super algebra, $\g$ is a Batalin-Vilkovisky
algebra with respect to the operator $\Delta=\bE$.
\end{proposition}
\begin{proof}

We know that $\bE^2=0$, so it suffices to check $\bE (\G \GH\GK)$ is of the
required form. This follows essentially because each term of $\bE$ can only
affect at most two of $\{\G,\GH,\GK\}$.
\end{proof}

This provides an alternate proof that $(-1)^{x} [\G,\GH]$ is a graded
Lie bracket
because of the following proposition. (See \cite{getzler},
\cite{chas}.)

\begin{proposition}\label{bvgerst}
Any BV algebra is also a Gerstenhaber algebra (in the super sense), by
defining the Lie bracket to satisfy:
$$(-1)^{|v|}[v,w] = \Delta(vw) - \Delta(v) w - (-1)^{|v|} v\Delta (w).$$
\end{proposition}

\section{Coproduct and cobracket}

In addition to the product structure $\mu\colon \g\otimes \g\to\g$ on the
graded vector space $\g$ of graphs, there is also a coproduct
structure
$\Delta\colon\g\to\g\otimes\g$.  To describe this, note that the algebra
structure on $\g$ induces an algebra structure on $\g\otimes\g$, by
$$(\G\otimes \GH) (\GK\otimes \GL)=(-1)^{yw}\G\GL\otimes \GH\GK.$$ 
For a connected
graph
$\G$, we define
$\Delta(\G)=\G\otimes 1+1\otimes\G$, where $1$ denotes the empty graph.
$\Delta$ is extended multiplicatively to all of $\g$ by
$\Delta(\G\cdot \GH)=\Delta(\G)\cdot\Delta(\GH)$.
For example, if $\G$ and $\GH$ are connected, then
$$\Delta(\G\cdot \GH)=\G\cdot \GH\otimes 1+1\otimes \G\cdot \GH+\G\otimes
\GH+(-1)^{xy} \GH\otimes \G .$$

As usual, it is convenient to express
things diagramatically when defining ``co-"objects.  If we extend a
map
$d\colon\g\to\g$ to
$d\colon\g\otimes\g\to\g\otimes\g$ by the formula
$d(\G\otimes \GH)=d(\G)\otimes
\GH + (-1)^{x}\G\otimes d(\GH)$, the fact that $d$ is a derivation can be
said efficiently as  $d\mu=\mu d$, i.e. the following diagram
commutes:

$$
\begin{matrix}
\g\otimes \g&{\buildrel \mu\over\to}&\g\\
d\downarrow  &&\downarrow d\\
\g\otimes\g&{\buildrel\mu\over\to}&\g
\end{matrix}$$

To define a {\it coderivation}, we simply reverse the arrows and
replace multiplication by
comultiplication: a map
$\delta$ is said to be a   co-derivation if $\Delta\delta =\delta
\Delta$, where $\delta$ is extended to $\g\otimes\g$
as before by $\delta(\G\otimes \GH)=\delta(\G)\otimes
\GH + (-1)^{x}\G\otimes \delta(\GH)$:
$$\begin{matrix}
\g\otimes \g&{\buildrel \Delta\over\leftarrow}&\g\\
\delta\uparrow  &&\uparrow\delta\\
\g\otimes\g&{\buildrel\Delta\over\leftarrow}&\g
\end{matrix}$$

\begin{proposition} \label{coderivation} $\bG$ is a coderivation.
\end{proposition}

\begin{proof}

    If $\G$ is connected, then all terms of
$\bG(\G)=\sum_e\G_e$ are connected, so that
\begin{align*}
\Delta(\bG\G )&=\Delta(\sum_e X_e)\\
&=\sum_e (\G_e\otimes 1 + 1\otimes  \G_e)\\
&=\bG\G\otimes 1 + 1\otimes \bG\G\\
&=\bG(\G\otimes 1 + 1\otimes \G)\\
&=\bG\Delta (\G)
\end{align*}
The case when $\G$ is not connected can be handled as follows. Since $\bG$ is a
derivation of $\g$, it is a derivation of $\g\otimes\g$. That is, for 
any $a,b,c$ and $d$ in $\g$,
$$\bG (a\otimes b)(c\otimes d) = (\bG (a\otimes b) ) (c\otimes d)
+ (-1)^{|a|+|b|}(a\otimes b)\bG (c\otimes d).$$
Apply this formula to $\bG \Delta(\G \GH) = \bG (\Delta(\G)\Delta(\GH))$,
using that you inductively know $\bG\Delta(\G) = \Delta\bG(\G)$ and 
$\bG\Delta(\GH) = \Delta\bG(\GH)$.
\end{proof}

\medskip

\begin{figure}[ht]
\begin{center}
\includegraphics{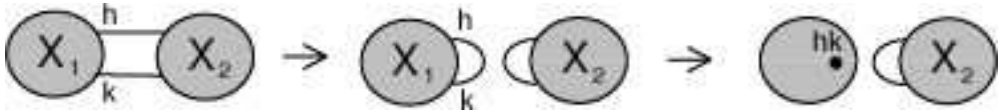}
\caption{Separating pair of half-edges}\label{separating}
\end{center}
\end{figure}

When we try the same computation with $\bE$, we run into problems
because the terms $\G_{hk}$ in $\bE(\G)$ may not be connected, even
when
$\G$ is connected, and the first line of the   calculation in
the proof of Proposition \ref{coderivation} is not valid.
For example, in  Figure \ref{separating}, contracting the half-edges 
$h$ and $k$ separates
the graph into two pieces.

By analogy with our definition of  the bracket, we define the {\it cobracket}
$\theta\colon
\g\to
\g\otimes
\g$ to measure how far
$\bE$ is from being a coderivation, i.e.

\begin{definition}  For any   graph $\G$, the {\it cobracket}
$\theta(\G)$ is defined to be
$$\theta(\G)= \Delta (\bE\G) -\bE \Delta(\G).$$
\end{definition}

\begin{figure}[ht]\begin{center}
\includegraphics{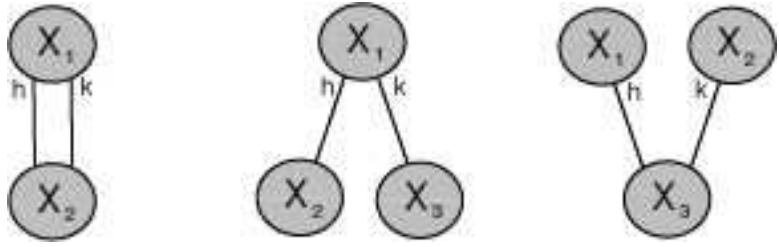}
\caption{All configurations of a separating pair}\label{sephk}
\end{center}
\end{figure}

In terms of graphs, the cobracket has the following interpretation.
If $\G$ is connected, we say that a pair $\{h,k\}$ of half-edges {\it
separates $\G$} if $\G_{hk}$ is not connected.  If $\{h,k\}$ 
separates a connected graph $\G$, then $h$
and $k$ must be in one of the
  configurations depicted in Figure \ref{sephk}, where the graphs 
$\G_i$ are connected.
We compute
\begin{align*} \Delta(\bE\G)&= \sum_{\{h,k\}|k\neq\bar h}\Delta(\G_{hk})\\
&=\sum_{\G_{hk} connected}(\G_{hk}\otimes 1 + 1\otimes \G_{hk}) +
\sum_{\G_{hk}=A_{hk}\cdot B_{hk} }\Delta(A_{hk}\cdot B_{hk})
\end{align*}

In the second summand, note that the graphs $A_{hk}$ and $B_{hk}$ are 
connected. If $A_{hk}$
has
$a$ vertices and
$B_{hk}$ has
$b$ vertices, we have
\begin{align*}
\Delta( A_{hk}\cdot B_{hk})=\Delta( B_{hk}\cdot A_{hk})&=
  A_{hk}\cdot B_{hk}\otimes1+1\otimes A_{hk}\cdot B_{hk}\\
&\qquad\qquad
+ A_{hk}\otimes B_{hk}+(-1)^{ab} B_{hk}\otimes A_{hk}
\end{align*}

On the other hand,
\begin{align*}
\bE\Delta(\G)&=\bE(\G\otimes 1+ 1\otimes\G)\\&=  \bE\G\otimes 1
+ 1\otimes \bE\G\\
&=\sum_{\{h,k\}|k\neq\bar h} (\G_{hk}\otimes 1 + 1\otimes\G_{hk})
\end{align*}
Thus the difference $\theta(\G)$ between $\bE\Delta(\G)$ and 
$\Delta(\bE\G)$ is the sum, over all pairs $\{h,k\}$ of half-edges
which separate $\G$, of
$ A_{hk}\otimes B_{hk}+(-1)^{ab} B_{hk}\otimes A_{hk}$.
We can simplify the notation by writing this
in the symmetric algebra, as

$$\theta(\G)=\sum_{\{h,k\}|\G_{hk}= A_{hk}\cdot B_{hk}} A_{hk}\odot B_{hk}.$$

\begin{figure}[ht]
\begin{center}
\includegraphics{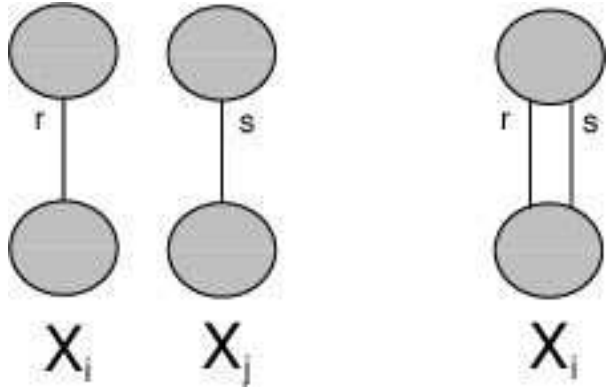}
\caption{Separating pairs in disconnected graphs}\label{Xij}
\end{center}
\end{figure}

     If $\G$ is not connected, the formula is more complicated.
Specifically, if $\G=\G_1\cdot \G_2\ldots\cdot \G_k$, with all $\G_i$
connected, we need to
consider  separating pairs $\{r,s\}$ in $\G_i$ and pairs $\{r,s\}$
with $r$ separating in $\G_i$ and $s$ separating in $\G_j$ (see 
Figure \ref{Xij}).
  Given such
a pair, write
$(\G_i)_{rs}=A\cdot B$ or $(\G_i\cdot \G_j)_{rs}=A\cdot B$, and list all
ordered partitions $(I,J)$ of $\{1,\ldots,\hat i,\ldots,\hat
j,\ldots, k\}$ into
two subsets (either or both of which may be empty). Then
$$\theta(\G_1\cdot \G_2\ldots\cdot \G_k)=\sum_{\{r,s\}}\sum_{(I,J)
}(-1)^\kappa \G_I\cdot A\odot B\cdot \G_J,$$
where $\G_I$ is the product of the $\G_i$ with $i\in I$, $\G_J$ is the
product of the $\G_j$ with $j\in J$, and $\kappa$ is the Koszul sign.
For example, if $k=2$, with  $\G_1=X$ and $\G_2=Y$, then
\begin{align*}
\theta(X\cdot Y)&=\sum_{r \in X, s\in Y, (X\cdot Y)_{rs}=A\cdot B}
A\odot B \\
&\qquad+ \sum_{r,s\in \G, \G_{rs}=A\cdot B} A\odot B\cdot Y 
+(-1)^{ab} B\odot A\cdot
Y\\ &\quad\qquad+\sum_{r,s\in Y, Y_{rs}=A\cdot B} X\cdot A\odot B +
(-1)^{ab} X\cdot B\odot A\\
\end{align*}
If we assume our graphs have
no separating edges, then the first summand above vanishes, and  the 
formula takes the following more elegant
form:  $$\theta(X\cdot Y)=\theta(X)\Delta(Y) +
(-1)^{xy}\theta(Y)\Delta(X),$$
where $X$ has $x$ vertices and $Y$ has $y$ vertices.  In fact this 
formula holds for arbitrary graphs $\G$ and $\GH$, not necessarily 
connected.
Here we have used the fact that $\g$ is graded cocommutative, hence
the coproduct induces a map $\g\to\g\odot\g$.

\begin{proposition}\label{abc} $$\theta^\phi\colon\g[-1]\to \Lambda^2\g[-1]$$ is a graded
Lie cobracket.
\end{proposition}

For the purposes of the proof it is easier to dualize. Since $\g$ is
the direct sum $\oplus\, \g_{v,e}$ of vector spaces spanned by graphs
with $v$ vertices and $e$ edges, and each $\g_{v,e}$ is finite
dimensional, from our
earlier remarks it follows that an operation $\theta\colon \g\to S^2 \g$ induces a
a graded Lie cobracket $\theta^\phi$ if and
only if
$\theta^\dagger=\oplus 
\theta_{v,e}^*$ satisfies the graded Jacobi identity on $\g^\dagger =\oplus\,
\g_{v,e}^*.$
Suppose $\G\in\g$. Let
$<\G,\cdot >$ denote the unique functional such that, for any graph $\GH$
$$<\G,\GH> =\begin{cases} |Aut(\G)|, & \text{if $\G\cong \GH$}\\
                    0 & \text{otherwise}\end{cases}$$
We denote the resulting isomorphism by $A\colon\g\to \oplus\,\g^\dagger$.

\begin{figure}[ht]
\begin{center}
\includegraphics{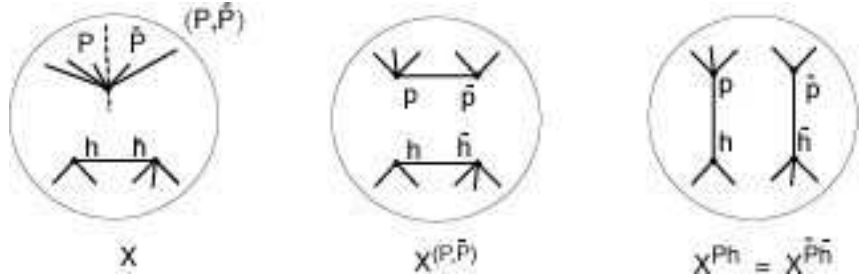}
\caption{Decontracting a partition and a half-edge}\label{decontract}
\end{center}
\end{figure}

The boundary operator $\bE$
will be replaced in the dual setting by a coboundary operator
$\delta_H\colon \g_v\to \g_{v+1}$,
defined as follows:

Let $\G$ be a graph, and let $(P,\bar P)$ be a   partition
of the edges incident to a vertex $v$
of
$\G$.   Expand the vertex $v$ to
obtain a new graph
$\G^{(P,\bar P)}$, with a new edge separating $P$ from $\bar P$.
This new edge is the union of two
half-edges, which we name $p$ and $\bar p$ to reflect the original
edges of $\G$ to which they are now incident.
The orientation on
$\G^{(P,\bar P)}$ is chosen so that collapsing the new edge gives back the
original orientation on $\G$. Given a half-edge $h$ in $\G$, we can 
now form the graph
   $\G^{(P,\bar P)}\langle ph \rangle$ for the graph obtained from 
$\G^{(P,\bar P)}$ by cutting and
pasting together the two half-edges $p$
and $h$.  If both $P$ and $\bar P$ have at least two elements, we 
denote this new graph by $\G^{Ph}$; otherwise
we set $\G^{Ph}=0$ (see Figure \ref{decontract}).

\begin{lemma}
As oriented graphs, $\G^{Ph}=\G^{\bar P\bar h}.$
\end{lemma}

\begin{definition} $\delta_H: \g_v\to \g_{v+1}$ is defined by
$$\delta_H(\G) = {1\over 2}\sum_{P, h} \G^{Ph},$$
where  $h$ runs over all
half edges of $\G$, and $P$ over all subsets of the edges at all vertices.
\end{definition}

  The factor of ${1\over 2}$ is there to account for the fact that, 
since $\G^{Ph}=\G^{\bar P\bar h}$, we have counted each graph in the 
coboundary
twice.

\begin{proposition}
The following diagrams commute.
\begin{equation*}
\begin{CD} \g^\dagger@>>{\partial_H^\dagger}>\g^\dagger\\
@AA{A}A @AA{A}A\\
\g@>>{\delta_H}>\g
\end{CD}
\hspace{5em}
\begin{CD} \g^\dagger\otimes\g^\dagger@>>{\Delta^\dagger}>\g^\dagger\\
@AA{A\otimes A}A @AA{A}A\\
\g\otimes\g@>>{\mu}>\g
\end{CD}
\end{equation*}
\end{proposition}
\begin{proof}

We first show that the left-hand diagram commutes.

Start in the lower-left with a graph $\G$.
Trace through
the diagram in both ways, evaluating in the upper-right hand corner on
the graph $\GH$.
Let $A_+$ denote the set $\{ (P,h): \G^{Ph} = \GH\}$.
Let $A_-$ denote the set $\{ (P,h): \G^{Ph} =  -\GH\}$.
Going right and then up in our diagram, we get
$$<\sum \G^{Ph}, \GH> = |Aut(\G^{Ph})|
(|A_+|-|A_-|)=|Aut(\GH)|(|A_+|-|A_-|).$$
Let $B_+$ denote the set $\{h,k: \GH_{hk} =\G \}.$
Let $B_-$ denote the set $\{h,k: \GH_{hk} = -\G \}.$
Going up and then right in the diagram:
$$<\G, \bE(\GH)> = |Aut(\G)|(|B_+|-|B_-|).$$

Thus it suffices to show that $|Aut(\GH)||A_\pm| = |Aut(\G)||B_\pm|$.
Now
$Aut(\G)$ acts transitively on each of $A_\pm$ and $Aut(\GH)$ acts
transitively on each of $B_\pm$, since we may assume that neither $\G$
nor $\GH$ have
orientation reversing automorphisms.
    Thus
$|Aut(\G)|/|stab(a_\pm)| = |A_\pm|$ and $|Aut(\GH)|/|stab(b_\pm)| =|B_\pm|$,
where
$stab(a_\pm)$,
$stab(b_\pm)$ denote the stabilizers of the elements $a_\pm\in A_\pm$ and
$b_\pm\in B_\pm$ respectively. It is thus sufficient to show that
$|stab(a_\pm)|=|stab(b_\pm)|$. To see this, suppose that $a_{\pm} =
(P,h)$.
    Notice that every automorphism of
$\G$ which fixes
$P$ and $h$ defines an automorphism of $\GH=\pm \G^{Ph}$.
Similarly every automorphism of $\GH$ which fixes $\{h,k\}$ extends to an
automorphism of
    $\G = \pm \GH_{hk}$. This gives us inverse maps between $stab(a_\pm)$ and
$stab(b_\pm)$.

Now we turn to the second commutative diagram. Start in the lower left
with a tensor $$\G_1^{n_1}\cdots \G_p^{n_p}\otimes
\G_1^{m_1}\cdots \G_p^{m_p},$$ where each $\G_i$ is a connected graph 
and $m_i,n_i\geq 0$.
To
establish commutativity, it suffices to evaluate in the upper right on
the graph
$$\G_1^{m_1+n_1}\cdots \G_k^{m_k+n_k},$$  since evaluating on other
monomials is zero in both directions. We may also assume that either
$\G_i$ has an even number of vertices, or that $m_i+n_i =1$, since 
otherwise the
oriented graph itself is zero. Therefore we may suppose that the 
tensor in the lower-left
corner is of the form
$$\G_1^{n_1}\cdots \G_k^{n_k}\GH_1\cdots \GH_s\otimes
\G_1^{m_1}\cdots \G_k^{m_k}\GH_{s+1}\cdots \GH_t$$
where $\G_i$ are distinct even graphs, and the $\GH_i$ are distinct odd
graphs. We evaluate in the upper right on the monomial
$\G_1^{m_1+n_1}\cdots \G_k^{m_k+n_k}\GH_1\cdots \GH_t$. To do this, when we go
up and then right, we must calculate
\begin{align*}
&\Delta(\G_1^{m_1+n_1}\cdots \G_k^{m_k+n_k}\GH_1\cdots \GH_t) =\ldots +\\
& \  \binom{m_1+n_1}{m_1}\cdots\binom{m_k+n_k}{m_k}\G_1^{n_1}\cdots
\G_k^{n_k}\GH_1\cdots \GH_s\otimes
\G_1^{m_1}\ldots \G_k^{m_k}\GH_{s+1}\cdots \GH_t\\
&\quad + \cdots\\
\end{align*}
Hence, going up and right we get 
\begin{align*}
&|Aut(\G_1^{n_1}\cdots \G_k^{n_k}\GH_1\cdots
\GH_s)|\cdot\\
&\qquad|Aut(\G_1^{m_1}\cdots \G_k^{m_k}\GH_{s+1}\cdots
\GH_t)|\cdot
\binom{m_1+n_1}{m_1}\cdots\binom{m_k+n_k}{m_k}.
\end{align*} 
 On 
the other hand,
going right and then up, we get simply $$|Aut(\G_1^{m_1+n_1}\cdots
\G_k^{m_k+n_k}\GH_1\cdots \GH_t)|.$$  Equality follows from the
fact that for any graph $\G$,
$|Aut(\G^{p+q})|=|Aut(\G^p)||Aut(\G^q)|\binom{p+q}{p}$; this is a 
consequence of
the formula $|Aut(\G^n)|= n!|Aut(\G)|^n$, coming from the fact that
$Aut(\G^n)$ is the semidirect product of the symmetric group 
$\Sigma_n$ and $Aut(\G)^n$.
\end{proof}

\medskip
\begin{proof}[Proposition \ref{abc}]

Note that by definition, $\theta = \Delta\partial_H-\partial_H\Delta
= \Delta\partial_H
-(\bE\otimes id + \tau(\bE\otimes id)\tau)\Delta$,
and therefore
$\theta^\dagger 
=\partial^\dagger_H\Delta^\dagger-\Delta^\dagger(\bE^\dagger\otimes 
id + \tau(\bE^\dagger\otimes id)\tau)$.
Hence $\hat{\theta}=A^{-1}\theta^\dagger(A\otimes A) = \delta_H\mu -
\mu(\delta_H\otimes id + \tau(\delta_H\otimes id)\tau)$ is the 
deviation from $\delta_H$ being a
derivation. Now by the first commutative diagram in the previous
proposition, $\delta_H^2 =0$. Also, $\delta_H$ satisfies the BV
axiom, since it is of degree $1$ and involves summing over pairs of
sub-objects of graphs (just like
$\partial_H $). Therefore by Proposition~\ref{bvgerst}, this deviation satisfies the graded
Jacobi identity:
$$\hat{\theta}\circ(\hat{\theta}\otimes id)\circ
(id+\sigma+\sigma^2) = 0.$$
This implies that $\theta^\dagger$ satisfies the graded Jacobi identity by the
following commutative diagram:
\begin{equation*}
\begin{CD}
\g^\dagger\otimes\g^\dagger\otimes\g^\dagger@>>{id +\sigma+\sigma^2}>
\g^\dagger\otimes\g^\dagger\otimes\g^\dagger@>>{\theta^\dagger\otimes
id}>\g^\dagger\otimes\g^\dagger@>>{\theta^\dagger}>\g^\dagger\\
@AA{A\otimes A\otimes A}A @AA{A\otimes A\otimes A}A @AA{A\otimes A}A
@AA{A}A\\
\g\otimes\g\otimes\g@>>{id + \sigma+\sigma^2}> \g\otimes\g\otimes\g
@>>{\hat{\theta}\otimes id}>
\g\otimes\g@>>{\hat{\theta}}>\g
\end{CD}
\end{equation*}
Now dualizing, $\theta$ satisfies graded co-Jacobi, which implies that
$\theta^\phi$ is a graded Lie cobracket.
\end{proof}

\section{Compatibility of bracket and cobracket}

Recall that if the
bracket $b=[\cdot,\cdot]$ and cobracket
$\theta$ satisfy the compatibility relation
$$\theta[v,w]+[\theta(v),w]+(-1)^{|\theta||v|}[v,\theta(w)]=0,$$
then the shifted complex $\g[-1]$ is a graded Lie bialgebra.

\begin{figure}[ht]
\begin{center}
\includegraphics{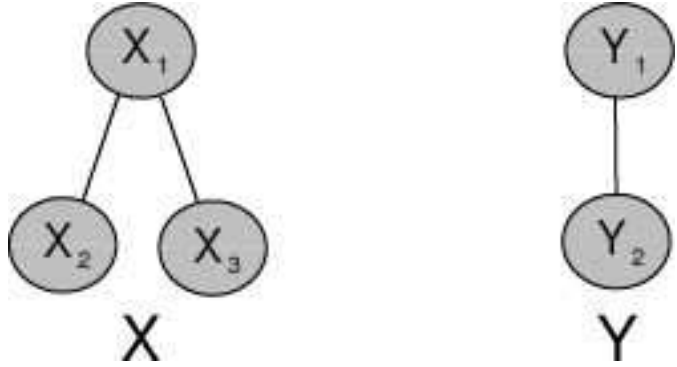}
\caption{Bracket and cobracket not compatible}\label{notcompat}
\end{center}
\end{figure}

In our case, the cobracket $\theta$ is of degree -1  (though $\theta^\phi$ is of degree
$0$).  The graph complex
$\g$ does not satisfy the compatibility relation,  as one can verify by computing using the 
graphs $\G$
and $\GH$ depicted in Figure \ref{notcompat}
The problem in this example is that both graphs $\G$ and $\GH$ have
separating edges, so that some terms in the bracket $[ \G,\GH]$
are not connected.  To remedy this, we consider the subspace
$\h\subset
\g$ generated by connected graphs with no separating edges.
In the literature these are often called {\it one-particle irreducible} graphs.
Note that $\h$ is a
subcomplex of $\g$ with respect to the boundary operator $\bG$,
though it is not with
respect to $\bE$.

The following two lemmas are easily verified:

\begin{lemma} If $\G$ and $\GH$ are connected graphs, and $\G$ has no
separating edges, then each term $(\G\cdot~\GH)_{hk}$ in the bracket
$[\G,\GH]$ is
connected.  If in addition $\GH$ has no separating edges, then each
term has no separating edges.
\end{lemma}

\begin{lemma} If $\G$ has no separating edges, and $\G_{hk}$ has two
components, then each component has no separating edges.
\end{lemma}

\begin{theorem} The subcomplex $\h\subset \g$ spanned by one-particle 
irreducible graphs is a graded Lie bialgebra with
respect to
$b^\phi=[\cdot,\cdot]^\phi$ and $\theta^\phi$.
\end{theorem}

\begin{proof}

    The Lemmas show that bracket and cobracket restrict
to operations on $\h$.  To prove the compatibility relation, we will show
that all terms in the sum
$\theta[\G,\GH]+[\theta(\G),\GH]+(-1)^{x}[\G,\theta(\GH)]$ cancel in 
pairs.  As usual,
$x$ denotes the number of vertices of $X$.

Because $\h$ has a basis consisting of connected graphs,
multiplication actually induces an isomporphism
$\mu:\h\odot\h\to\h^2$. Composing $\mu$ with $\theta$ simplifies the
expression for the cobracket to:
    $$\mu \theta (\G)=\sum_{\{h,k\} separating } \G_{hk}.$$
We also have   $\mu[\G\odot \GH,\GK] = [\G\GH,\GK]$ and $\mu[\G,\GH\odot
\GK]=[\G,\GH\GK]$, since the bracket satisfies the
Gerstenhaber axiom, so that  showing that
$\theta[\G,\GH]+[\theta(\G),\GH]+(-1)^{x}[\G,\theta(\GH)]$ vanishes is
equivalent to showing that
$$\mu\theta[\G,\GH] +[\mu\theta(\G),\GH] + (-1)^{x}[\G,\mu\theta(\GH)] = 0.$$

We have
$$\mu\theta[\G,\GH]=\mu\theta(\sum_{h\in \G, k\in \GH}(\G\cdot 
\GH)_{hk})=\sum_{h\in
\G, k\in \GH}\sum_{r,s}((\G\cdot \GH)_{hk})_{rs},$$
where $r,s$ are separating half edges in $(\G\cdot \GH)_{hk}$.

For each term $((\G\cdot \GH)_{hk})_{rs}$ of $\mu\theta[\G,\GH]$ there are
several cases.

{\bf Case 1}.  $r,s\subset \G,  r,s\neq \bar h$.

If $\{r,s\}$ is not a separating pair in $\G$, then the edge $e(k)$ 
containing $k$ must be a separating edge in $Y$, contradicting our
assumptions.  Thus  $(\G\cdot \GH)\langle hk \rangle$ are arranged as 
in Figure \ref{case1},
where $e(r)$ contains $r$  and $e(s)$ contains $s$.

\begin{figure}[ht]\begin{center}
\includegraphics{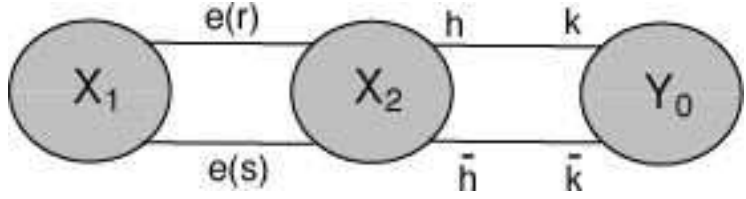}
\caption{$(\G\cdot\GH)\langle hk\rangle$ in Case 1 ($r,s\subset \G, 
r,s\neq \bar
h$)}\label{case1}
\end{center}\end{figure}

Then $\G_{rs}$ is a term of $\mu\theta (\G)$, and
$(\G_{rs}\GH)_{hk}=((\G\GH)_{rs})_{hk}$ is a term of $[\mu\theta(\G),\GH]$,
which cancels with $((\G\GH)_{hk})_{rs}$ by Lemma \ref{orient}.

{\bf Case 2}.  $r,s \subset \GH, r,s\neq \bar k$.

This is similar to the last case (see Figure \ref{case2}).
\begin{figure}[ht]\begin{center}
\includegraphics{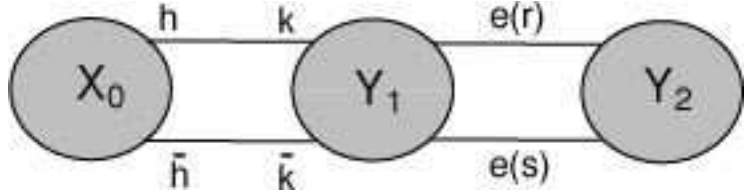}
\caption{$(\G\cdot\GH)\langle hk\rangle$ in Case 2 ($r,s \subset \GH, 
r,s\neq \bar
k$)}\label{case2}
\end{center}\end{figure}
$\GH_{rs}$ is a term of $\mu\theta(\GH)$ and $(\G\cdot \GH_{rs})_{hk} =
(-1)^{x} ((\G\GH)_{rs})_{hk}$ is a term of $(-1)^{x}[\G,\mu\theta(\GH)]$
which cancels with $((\G\GH)_{hk})_{rs}$.

{\bf Case 3}. $r\subset \G, s\subset \GH, r\neq \bar h, s\neq\bar k$

Here $(\G\cdot \GH)\langle hk \rangle$ must be in the configuration 
of Figure \ref{case3},
since neither
$X$ nor
$Y$ has a separating edge.

\begin{figure}[ht]\begin{center}
\includegraphics{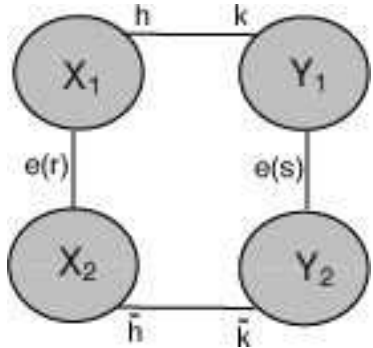}
\caption{$(\G\cdot\GH)\langle hk\rangle$ in Case 3 ($r\subset \G, 
s\subset \GH, r\neq \bar h,
s\neq\bar k$)}\label{case3}
\end{center}\end{figure}

But then $(\G\GH)_{rs}$ is a term of $[\G,\GH]$ and $((\G\GH)_{rs})_{hk}$ is a
term of $\mu\theta[\G,\GH]$ which cancels
with $((\G\GH)_{hk})_{rs}$.

{\bf Case 4}. $r =\overline{h}$

$(\G\cdot \GH)\langle hk \rangle$ must be as in Figure \ref{case4}.
\begin{figure}[ht]\begin{center}
\includegraphics{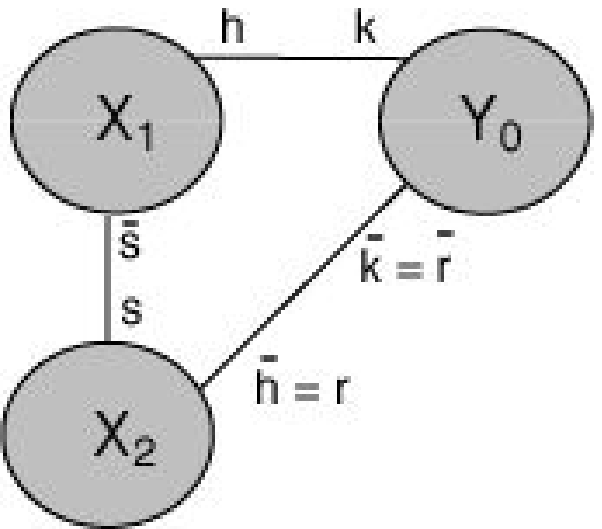}
\caption{$(\G\cdot\GH)\langle hk\rangle$ in Case 4 ($r 
=\overline{h}$)}\label{case4}
\end{center}\end{figure}

In this case $\G_{rs}$ is a term of $\mu\theta(\G)$, and so
$(\G_{rs})_{hk}$ is a term of $[\mu\theta(\G),\GH]$
which cancels.

{\bf Case 5}.  $r=\overline{k}$

This is similar to case $4$ (see Figure \ref{case6}.

\begin{figure}[ht]\begin{center}
\includegraphics{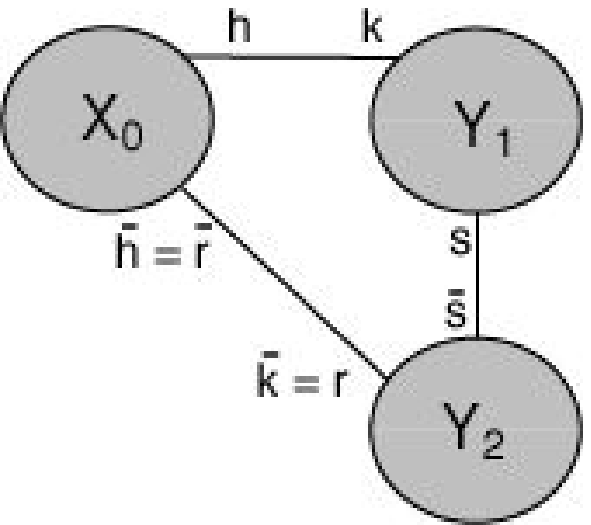}
\caption{$(\G\cdot\GH)\langle hk\rangle$ in Case 6 ($r 
=\overline{k}$)}\label{case6}
\end{center}\end{figure}

Here $(-1)^{x}((\G\GH)_{rs})_{hk} = (\G\cdot \GH_{rs})_{hk}$ is a term of
$(-1)^{x}[\G,\mu\theta(\GH)]$ which
cancels.

Thus we have shown that each term of $\mu\theta[\G,\GH]$ cancels either
with another term of $\mu\theta[\G,\GH]$ or with a term of 
$[\mu\theta(\G),\GH]$
or
$(-1)^h[\G,\mu\theta(\GH)]$. Also, each term of $[\mu\theta(\G),\GH]$  and
$(-1)^{x}[\G,\mu\theta(\GH)]$ does, in
fact, cancel with a term in $\mu\theta[\G,\GH]$.
\end{proof}

\noindent{\bf Remark:} 
The technique of this proof can be used to give an easy
argument that $\theta^2=0$ on
$\mathcal H$.
\medskip

The fact that the subcomplex $\mathcal H$ carries a bi-algebra
structure is valid in
the  general setting of cyclic operads.  For the associative and Lie
cases,
$\mathcal H$ is quasi-isomorphic to  $Prim(\g)$
(see\cite{exposition}), so that we have a bi-algebra structure on chain
complexes $\mathcal H$ which compute the cohomology of mapping class
groups (associative operad) and of
 groups of outer automorphism groups of free groups  (Lie operad).
In the commutative case, which we have focused on in this paper,
computer calculations
due to F. Gerlits
\cite{gerlits} show that the map $H_*(\h)\to
H_*(Prim(\g))$ is not surjective.  However, in \cite{cv3} we show that
the Lie bracket and cobracket described in this paper do induce a Lie
bi-algebra structure on a certain quotient complex $Prim(\g)/\mathcal C$
which is quasi-isomorphic to $Prim(\g)$, so in this case as well we have a
bi-algebra structure on a chain complex which computes graph homology.

\section{Homology}

In this section we consider $\g$ as a complex with boundary operator $\bG$.

\begin{proposition} The graph bracket descends to the level of homology.
\end{proposition}

\begin{proof}

This follows from the identity
$$\bG[\G,\GH]=[\bG \G, \GH]+(-1)^{x}[\G,\bG \GH],$$
Which can be derived by expanding the identity
$$(\bE\bG+\bE\bG)(\G \GH) = 0.$$
\end{proof}

\begin{proposition}
The graph cobracket descends to the level of homology.
\end{proposition}
\begin{proof}

The proof of the previous lemma can be dualized to yield
a proof of the present lemma by considering the equation
$$\Delta ((\bE\bG+\bG\bE)\G)=0.$$
In the course of expanding this out one must use the fact that
$\bG$ is a coderivation and that
$\bE\otimes id + id\otimes \bE$ anti-commutes with
$\bG\otimes id + id\otimes \bG$, where one must
as always respect the Koszul rule of signs. \end{proof}

It turns out that the bracket can be interpreted in terms of the Lie
algebra $c_n$ mentioned in the introduction.
Kontsevich constructs an isomorphism between
$\mathfrak{sp}(2n)$-invariants in $\Lambda^*c_n$ and elements of $\g$.
 More specifically, there is a map
$$\phi_n: (\Lambda c_n)^{\mathfrak{sp}(2n)}\to \g.$$
See \cite{exposition}.

\begin{proposition}
$\phi_n$ is a Lie algebra homomorphism. That is, it maps the Schouten
bracket to the graph bracket.
\end{proposition}

\begin{proof}

We recall that an $\mathfrak{sp}(2n)$ invariant tensor is
associated to a graph by the following procedure. Each vertex
of the graph represents a tensor factor, in the order given by
the vertex labelling. For
each edge we put a $p_i$ at the tail of the arrow and a $q_i$
at the head or we put the $q_i$ at the tail and the $p_i$ at
the head, incurring a minus sign as a result. We sum over all
possible choices, each choice is called a ``state.'' Passing to
the wedge product yields an $\mathfrak{sp}(2n)$ -invariant. The
Schouten bracket
involves first choosing two tensor factors to bracket, which means 
picking a vertex from each of the two graphs.
Then one takes the Poisson bracket of the monomials
at each vertex. This can be thought of as deleting a $p_i$
from one and a $q_i$ from the other, and then multiplying
the monomials together. One can view the result as breaking
the edges with the $p_i$ and $q_i$ into half-edges, gluing
them together and contracting, and also gluing together the
resulting dangling edges. Summing over all possible states, we
see that this contribution to the Schouten bracket is given
by contracting the two given half-edges, which is the definition
of the bracket.
\end{proof}

It is not difficult to show that the Schouten bracket is always trivial on
the homology level: one can think of it as the deviation of the 
Chevalley-Eilenberg differential from being a derivation with respect
to the wedge product. Therefore, the preceding proposition
suggests that
the bracket is trivial on the homology level.
This is not quite true, however, since $\phi_n$ is not a chain map!
(This was an oversight in Kontsevich's argument. In \cite{exposition}
we show how to repair this oversight.)
 However, there is a straightforward
proof of the homological triviality of the  bracket, due to S.
Mahajan.

\begin{proposition}
The bracket $[\cdot,\cdot]$ is trivial on the homology level.
\end{proposition}
\begin{proof}

Define a multiplication $\mu_1:\g\otimes\g\to\g$ by $$\mu_1(X\otimes Y) =
\sum_{x\in X,y\in Y} (X\cdot Y)<xy>$$
(See Figure \ref{halfcontract}.) 
Then an easy argument shows that $[X,Y]= \partial_E\circ \mu_1-\mu_1\circ\partial_E$.
Therefore, if $X$ and $Y$ are both cycles, so is $X\cdot Y$, and
$[X,Y]=\partial_E \mu_1(X\otimes Y)$. \end{proof}

A similar argument shows

\begin{proposition}
The cobracket $\theta$ is trivial on the homology level.
\end{proposition}
\begin{proof}
This follows by a very similar argument to the previous proposition.
Define $\Delta_1:\g\to\g\otimes\g$ in the same way as $\theta$ only don't contract
an edge. Then $\theta = \partial_E\circ \Delta_1-\Delta_1\circ\partial_E$.
\end{proof}

Finally, we show that $\partial_H$ is zero at the homology level as well.

\begin{proposition}
$\partial_H$ is zero on homology.
\end{proposition}
\begin{proof}

Define a map $\alpha\colon \g\to \g$ as
$$\alpha (X) =\frac{1}{2} \sum_{x,y, x\neq\bar{y}} X<xy>.$$
It is straightforward that $\partial_E\alpha - \alpha
\partial_E=\partial_H$. 
\end{proof}

\end{document}